\documentclass[11pt]{article}
\usepackage{epsfig}
\usepackage{amsmath}
\usepackage{amsfonts}
\usepackage{amssymb}
\usepackage{color}

\title{The pluri-fine topology is locally connected}

\author{Said El Marzguioui and Jan Wiegerinck\\}

\newcounter{bean}

\newtheorem{theorem}{Theorem}[section]
\newtheorem{corollary}[theorem]{Corollary}
\newtheorem{lemma}[theorem]{Lemma}

\newcommand{\qed}{\hspace*{\fill}$\square$}

\DeclareMathOperator*{\f-limsup}{f-\,limsup}
\begin{document}

\maketitle \footnote{2000 Mathematics Subject Classification 32U15, 31C40}

\begin{abstract}
We prove that the pluri-fine topology on any open set $\Omega$ in $\mathbb{C}^{n}$ is locally
connected. This answers a question by Fuglede in [4]. See also Bedford [6].
\end{abstract}
\begin{center}
\section{Introduction}
\end{center}
The pluri-fine topology on an open set $\Omega$ in $\mathbb{C}^{n}$ is the coarsest topology on
$\Omega$ making all the plurisubharmonic functions on $ \Omega$ continuous. Almost all the results
concerning the classical fine topology, introduced by H. Cartan in 1940, remain valid, and even
with the same proofs. For example, it is obvious that the pluri-fine topology is Hausdorff,
completely regular. It was observed by Bedford and Taylor [5], that it is Baire (i.e. a countable
intersection of pluri-finely open and pluri-finely dense sets, is pluri-finely dense), and it has
the quasi-Lindel\"{o}f property (i.e. every family of pluri-finely open sets contains a countable
subfamily whose union differs from that of the whole family at most by some pluripolar set). These
results were discovered in the classical case by Doob [7] in 1966. See Bedford and Taylor [5] and Klimek [11].\\
In 1969 Fuglede [1] showed that the classical fine topology is locally connected, and that any
usual domain (connected open set) is also a fine domain. This was for him the starting point for
developing an interesting theory of finely harmonic functions, defined on finely open set (see e.g.
\cite{Fug72}). Few years later, Fuglede and others have developed the theory of finely holomorphic
functions on $ \mathbb{C}$ which recently has found a surprising application in
pluripotential theory. See Edlund and J\"{o}ricke, [13].\\
In [4] Fuglede made an attempt to introduce fine holomorphy in $ \mathbb{C}^{n}$. He compared three
possible fine topologies on $ \mathbb{C}^{n}$, the fine topology on $ \mathbb{R}^{2n}$, the
pluri-fine topology and the $n$-fold product topology induced by the fine topology on $
\mathbb{C}$. He makes it clear that the pluri-fine topology is the right one to use. Then he notes
that local connectivity needs to be established before fine holomorphy can be developed at all.
\\The proof given by Fuglede in [1] of the local connectivity of the fine topology in $ \mathbb{R}^{n}$ was strongly
based on the theory of balayage of measures, especially the balayage of the unit Dirac measure.
Unfortunately, no such theory exists in pluripotential theory and it seems to be impossible to
develop, because the strong subadditivity of the relative extremal plurisubharmonic function fails
to hold, as was proved by Thorbi\"{o}rnson [8]. Moreover, unlike the situation in classical
potential theory, the notions of thinness and pluripolarity are not equivalent. This means that
plurithin sets can not be characterized in terms of capacity, which accounts for a big difference
between the pluri-fine and fine topology. However, using elementary properties of finely
subharmonic functions, that were found by Fuglede [2, 3], we give
a surprisingly simple proof of the local connectivity of the pluri-fine topology.\\
\begin{theorem}
The pluri-fine topology on an open set $\Omega$ in $\mathbb{C}^{n}$ is locally connected.
\end{theorem}
\begin{corollary}
Every pluri-finely connected component of a pluri-finely open set $ \Omega \subset \mathbb{C}^{n}$
is pluri-finely open. Moreover the set of these components is at most countable.
\end{corollary}
The proof of this corollary is the same as the proof of Fuglede in the classical case. It uses the
quasi-Lindel\"{o}f property. See Fuglede [1].\\ In order to make this note more complete, let us
mention the following two theorems due to Fuglede [1] in the classical case.
\begin{theorem}An open set $U \subseteq \mathbb{C}^{n}$ is pluri-finely connected if and only if $U$
is connected in the Euclidean topology.
\end{theorem}
\begin{theorem}Let $U \subseteq \mathbb{C}^{n}$ be a pluri-finely open and pluri-finely connected set. If $E$ is a
pluripolar set, then $U \backslash E$ is pluri-finely connected.
\end{theorem}
Thanks to the fact that the fine topology on $\mathbb{R}^{2n}$ is finer than the pluri-fine
topology on $\mathbb{C}^{n}$, these two theorems are an immediate consequence of similar results
proved by Fuglede [1] in classical fine topology. See also Fuglede [2] for different proofs. \\
It should be mentioned that contrary to the classical fine topology, the set $U \backslash E$ of
Theorem 1.4 is not pluri-finely open in general.
\begin{center}
\section{Preliminaries}
\end{center}
We will need the following two results from [2, p. 100 and p. 87, respectively].
\begin{lemma} Suppose that $V \subseteq U \subseteq \mathbb{C}^{n}$ are finely open sets, and let $\psi$
(resp $\varphi$ ) be a finely subharmonic function on $U$ (resp $V$). Assume that:
$$
\f-limsup_{z \to x ,{z\in V}}\varphi (z) \leq \psi (x) \ for \ all \ x \in U\cap {\partial_{f}V}.
$$
 Then the following function
$\Psi$ is finely subharmonic in ${U}$:
$$
\Psi(z)= \left\{
\begin{array}{ll}
\max \{ \varphi (z), \psi (z)\} & \mbox{if $z \in V$},\\
\psi (z) & \mbox{if $z \in U \diagdown V$}.
\end{array} \right.
$$
\end{lemma}
Here $\f-limsup$ denotes the $\limsup$ with respect to the fine topology (i.e. $ \inf \limits_{O}
\sup \limits_{ z \in O} \psi (z)$ where $O$ ranges over the set of all fine open sets in $V$ which
contain $z$) and $\partial_{f}$ stands for the fine boundary.
\begin{theorem}Let $U$ be an open set in $\mathbb{C}^{n}$. Then a
function $\varphi$ :\ $U$ $\longrightarrow$ $\mathbb{R}$ is subharmonic iff $\varphi$ is finely
subharmonic and moreover locally bounded from above in the Euclidean topology.
\end{theorem}

\noindent {\sc \textbf{Remark}} It was proved by Fuglede [3] that in the case of the plane, i.e.
$n=1$, the local boundedness from above may be omitted in the theorem. He also gave examples which
prove that the condition can not be removed in higher dimensions. (see e.g.
\cite{Fug74}).\noindent\\
It is a fundamental result of H. Cartan (see e.g. \cite{Br71}), that pluri-fine neighbourhoods of
$z_{0}$ are precisely the sets of the form $ \mathbb{C}^{n}\backslash E$, where $E$ is pluri-thin
at $ z_{0}$ and $z_{0} \not\in E$. A subset $E$ of $\mathbb{C}^{n}$ is said to be pluri-thin at $
z_{0}$ if and only if either $ z_{0}$ is not a limit point of $E$ or there is $ r> 0$ and a
plurisubharmonic functions $\varphi$ on $B(0,r)$ such that
$$
\limsup_{z \to z_{0},z\in E\backslash \{z_{0}\}}\varphi(z)< \varphi(z_{0}).
$$
The following important result (see e.g. \cite{K91}) and its corollary assert, that being a
pluri-fine open is a local property.
\begin{theorem}
Finite intersections of sets of the form
$$
B_{\varphi,\Omega,c} = \{z\in \Omega :\varphi (z)> c\},
$$
where $\Omega \subseteq \mathbb{C}^{n}$ is open , $\varphi \in PSH(\Omega)$ and $c \in \mathbb{R}$,
constitute a base of the pluri-fine topology on $\mathbb{C}^{n}$.
\end{theorem}
\begin{corollary}If $ \Omega_{1} \subseteq \Omega_{2} \subseteq \mathbb{C}^{n}$ are open subsets,
then the pluri-fine topology on $ \Omega_{1}$ is the same as the topology on $ \Omega_{1} $ induced
by the pluri-fine topology on $ \Omega_{2}$.
\end{corollary}
\begin{center}
\section{Proof of theorem 1.1}
\end{center}
The following result is the key to our proof of theorem 1.1. it was stated by Bedford and Taylor in
[5, theorem 2.3].
\begin{lemma}
Sets of the form
$$
B_{ \Omega}^{ \varphi} = \{z\in \Omega :\varphi (z)> 0\},
$$
where $\Omega \subseteq \mathbb{C}^{n}$ is open, $\varphi \in PSH(\Omega)$, constitute a base of
the pluri-fine topology on $\mathbb{C}^{n}$.
\end{lemma}
Bedford and Taylor did not give the proof of this result. Instead they referred to [7, 10].
However, the results in [7, 10] assert only that, sets of the form $B_{ \Omega}^{ \varphi}$
constitute a subbases for the fine topology. Since we could not find a proof in the literature, we
give a
proof here.\\

\textit{Proof}. By theorem 2.3 $B_{ \Omega}^{ \varphi}$ is a pluri-finely open set. Let $U\subseteq
\mathbb{C}^{n}$ be a pluri-finely open set, and let $a \in U$. It is a basic theorem of H. Cartan
that the complement $E$ of $U$ is pluri-thin at $a$. We will prove that $U$ contains a pluri-finely
open neighbourhood of $a$ of the form stated in the lemma. This is trivial if $a$ belongs to the
Euclidean interior of $U$. Consider the case when a is an accumulation point of $E$. There exist
then $ \delta > 0$ and a plurisubharmonic function $ \varphi$ on $ B(a, \delta)$ such that
$$
\limsup_{z \to a,z\in E}\varphi(z)<\varphi(a).
$$
Without loss of generality we may suppose that $ \varphi ( E \cap B(a, \delta)) \leq 0 < \varphi
(a)$. Since $ \mathbb{C}^{n} \backslash (E \cap B(a, \delta)) = U \cap B(a, \delta) \cup
\mathbb{C}^{n} \backslash B(a, \delta)$, we get $ \{ z \in B(a, \delta) : \ \varphi (z) > 0 \}
\subset U \cap B(a,
\delta) \subset U$. Which proves the lemma. \qed  \\

Denote by $B= B(0,1)$ the open unit ball in $ \mathbb{C}^{n}$, and let $\varphi \in PSH(B(0,1))$
such that $0\leq\varphi\leq1$ on $ B(0,1).$
\begin{lemma}
Let $U = \{\varphi>0\} \cap B(0,1)$. If $U = V \cup W$, where $V$ and $W$ are non empty pluri-fine
open sets such that $V \cap W= \O$, then the following function is plurisubharmonic:
$$
\varphi_{V}(z)= \left\{
\begin{array}{ll}
\varphi (z) & \mbox{if $z \in B\backslash W$},\\
0 & \mbox{if $z \in W$}.
\end{array} \right.
$$
\end{lemma}
\textit{Proof}. Let $\zeta \in \partial B$ and denote by $L_{\zeta}$ the complex line through $0$
and $\zeta$. We will prove that the function,
$$
\varphi_{V,\zeta}(z)= \left\{
\begin{array}{ll}
\varphi_{\zeta} (z) & \mbox{if $z \in L_{\zeta}\cap B \backslash W_{\zeta}$},\\
0 & \mbox{if $z \in W_{\zeta}$},
\end{array} \right.
$$
is subharmonic on $L_{\zeta} \cap B$. Here $ \varphi_{\zeta}$ denotes the restriction of $\varphi$
to $ L_{\zeta} \cap B$, and $W_{\zeta}:= W \cap L_{\zeta}$ . The set $W_{\zeta}$ is a fine open
subset of
$L_{\zeta}\cap B$ which may be empty.\\
Note that $\varphi_{V,\zeta}$ is the restriction of $ \varphi_{V}$ to $L_{\zeta}\cap B$. Denote by
$\partial_{f} W_{\zeta}$ the fine boundary of $W_{\zeta}$ relative to $ L_{\zeta} \cap B$. We claim
that $\varphi_{\zeta} =0$ on $\partial_{f} W_{\zeta}$. To prove the claim observe first that
$\partial_{f}W_{\zeta} \subset
\partial_{f} \{\varphi_{\zeta}
> 0\}$ and $\partial_{f}\{ \varphi_{\zeta}> 0\}= \partial_{f}\{\varphi_{\zeta}=0\}$. Moreover, the
set $\{\varphi_{\zeta}=0\}$ is finely closed, which means that $\partial_{f}\{\varphi_{\zeta}=0\}$
is a subset of $\{\varphi_{\zeta}=0\}$, and hence $\partial_{f} W_{\zeta} \subset
\{\varphi_{\zeta}=0\}$.\\
Next, we can assume that $L_{\zeta}\cap B \backslash W_{\zeta}$ is nonempty, for otherwise $
\varphi_{V, \zeta} \equiv 0$ hence subharmonic. Using the claim and the fact that $ \varphi $ is a
non-negative upper-semicontinuous function, we get the following:\\
$$
\limsup_{z\to a ,z \in B\cap L_{\zeta} \backslash W_{\zeta}}\varphi_{\zeta}(z) \leq
\varphi_{\zeta}(a) = 0, \ \forall a \in
\partial_{f} W_{\zeta},
$$

and clearly,
$$
\f-limsup_{z \to a,z\in B\cap L_{\zeta}\backslash W_{\zeta}}\varphi_{\zeta}(z)\leq 0, \ \forall a
\in \partial_{f} W_{\zeta},
$$
because the ordinary $\limsup$ majorizes the $\f-limsup$.\\
In view of the claim, the definition of $\varphi_{V, \zeta}$ does not change if we replace
$W_{\zeta}$ by its fine closure $W_{\zeta}^{f}$. Since $ \partial_{f} (L_{\zeta} \cap B \backslash
W_{\zeta}) \cap B =\partial_{f} W_{\zeta} $, Lemma 2.1 applies and $\varphi_{V,\zeta}$ is therefore
finely
subharmonic, which is clearly bounded, and hence subharmonic by Theorem 2.2.\\
It is a well known result that a bounded function, which is subharmonic on each complex line where it is defined is a plurisubharmonic.
(see e.g. \cite{P45}) \qed  \\
\\
\textit{Proof of Theorem1.1}. Let $z_{0} \in \mathbb{C}^{n}$ and let $D$ be a pluri-fine open
neighborhood of $z_{0}$. By Lemma 3.1 there exist an open set $\Omega$ in $\mathbb{C}^{n}$ and a
plurisubharmonic function $\varphi \in PSH( \Omega)$, such that the set $\{z \in \Omega : \
\varphi(z)> 0\}$ is a pluri-fine open neighborhood of $z_{0}$ contained in $D$. In view of
Corollary 2.4 and the fact that the pluri-fine topology is biholomorphically invariant, there is no
loss of generality if we assume that $z_{0}=0$, $ \Omega$ is the unit ball $B(0,1)$ and that $0
\leq \varphi \leq 1$ on $B(0,1)$. To prove the Theorem we will find a pluri-fine open neighborhood
of $0$ which is pluri-finely connected and
contained in $B_{B(0,1)}^{ \varphi} := \{ z \in B(0,1) : \ \varphi(z) >0 \}$.\\
Denote by $\mathcal{F}$ the set of all pluri-finely open sets $V$ which contain $0$ and for which
there exists a non empty pluri-finely open set $W_{V}$ such that $B_{B(0,1)}^{ \varphi}= V \cup
W_{V}$ and $V \cap W_{V} = \O$. It follows from Lemma 3.2 that the function $\varphi_{V}$ is
plurisubharmonic for all $V \in \mathcal{F}$. Moreover, the family $( \varphi_{V})_{V \in
\mathcal{F}}$ is left directed and lower bounded by $0$. It is a classical result, (see e.g.
\cite{H69} Theorem 4.15), that the infimum $ \psi$ of such a family exists and is plurisubharmonic
on $B(0,1)$. Now we claim that the set $U = \{ z \in B_{B(0,1)}^{ \varphi} : \ \psi (z) > 0\}$ is
pluri-finely open and a
pluri-finely connected neighborhood of $0$.\\
To prove the claim observe first that for every $V$, $ \psi(0) = \varphi_{V}(0) = \varphi (0)> 0,$
which means that $U$ is a nonempty pluri-fine open neighborhood of $0$. Next, from the definition
of $\varphi_{V}$ in Lemma 3.2 we see that $\varphi_{V}(z)= \varphi (z) > 0$ on $V$ and $
\varphi_{V}(z)= 0 $ on $B(0,1) \backslash V$. Consequently, $U = \bigcap \limits_{V \in
\mathcal{F}} \{\varphi_{V}> 0\}= \bigcap \limits_{V \in \mathcal{F}}V$ and $ B_{B(0,1)}^{ \varphi}
= \bigcap \limits_{V \in \mathcal{F}}V \cup \bigcup \limits_{V \in \mathcal{F}}W_{V}$, where
$\bigcap \limits_{V \in \mathcal{F}}V \cap \bigcup \limits_{V \in \mathcal{F}}W_{V} =\O$. Therefore
$U$ is an element of $\mathcal{F}$. It is minimal in the sense that it can not be split into two
disjoint non-empty pluri-fine open sets, which proves the claim and the theorem.

\noindent \textsc{KdV Institute for Mathematics \\
University of Amsterdam \\
Plantage Muidergracht 24 \\
1018 TV Amsterdam \\
The Netherlands}\\
janwieg@science.uva.nl \\  smarzgui@science.uva.nl

\end{document}